\documentclass{ifacconf}
\pdfoutput=1

\usepackage[utf8]{inputenc}

\usepackage{graphicx}  
\usepackage{natbib}        
\usepackage{amsmath}
\usepackage{amssymb}
\usepackage{epstopdf}
\usepackage{mathtools}
\usepackage{letltxmacro}

\newcommand{\nref}[1]{(\ref{#1})}
\renewcommand{\bf}[1]{\boldsymbol{#1}}
\newcommand{\ncite}[1]{[\cite{#1}]}
\renewcommand{\endassum}{\hspace*{\fill} \qed}
\newtheorem{definition}{Definition}
\begin{document}

\begin{frontmatter}

\title{An integral Nash equilibrium control scheme for a class of multi-agent linear systems} 
\thanks[footnoteinfo]{This work was partially supported by the ERC under research project COSMOS (802348). E-mail addresses: \{s.krilasevic-1, s.grammatico\}@tudelft.nl.}

\author{Suad Krilašević} and
\author{Sergio Grammatico} 

\address{Delft Center for Systems and Control, TU Delft, The Netherlands}

\begin{abstract}                
We propose an integral Nash equilibrium seeking control (I-NESC) law which steers the multi-agent system composed of a special class of linear agents to the neighborhood of the Nash equilibrium in noncooperative strongly monotone games. First, we prove that there exist parameters of the integral controller such that the system converges to the Nash equilibrium in the full-information case, in other words, without the parameter estimation scheme used in extremum seeking algorithms. Then we prove that there exist parameters of the I-NESC such that the system converges to the neighborhood of the Nash equilibrium in the limited information case where parameter estimation is used. We provide a simulation example that demonstrates that smaller perturbation frequencies and amplitudes are needed to attain similar convergence speed as the existing state-of-the-art algorithm.

\end{abstract}


\end{frontmatter}

\section{Introduction}
Extremum seeking control is a class of data-driven, adaptive control techniques used in optimization problems where the cost is a function of the states of a dynamical system. The method is a \textit{zero-order} method which means it only uses the value of the cost function for optimization and no a priori knowledge of the cost function is needed, except for some basic assumptions.

The method was first proposed in \ncite{leblanc}. For many years there was no analytical proof of stability of extremum seeking control for general nonlinear systems until the paper \ncite{krstic2000stability}. This sparked renewed interest into further development of this type of control. Most of the research was based on the original paper by Krstić and Wang, e.g. \ncite{tan2006non}, \ncite{ghaffari2012multivariable}, etc. There were also methods based on different ideas, such as \ncite{durr2013lie}, where the authors proposed an extremum seeking scheme based on Lie algebra, which turned out to be equivalent to the Krstić-Wang scheme. Based on the parameter estimation scheme from \ncite{adetola2007parameter}, Guay and Dochain propose an extremum seeking scheme \ncite{guay2017proportional} which does not use singular perturbation and averaging theory. As a result, a faster convergence rate is obtained. This fact motivates further research on such type of extremum seeking controller.

Although very similar at a first glance, Nash equilibrium problems (NEP) are different from standard (distributed) optimization problems, as they are characterized by a number of selfish agents whose goal is to optimize their individual cost functions, each possibly dependent on the decision variables of other agents. In NEPs, the constraints on the decision set of each agent are independent of other agents, while in the generalized Nash equilibrium problems (GNEP), they share constraints. Recent interest in GNEPs is justified by the fact that many engineering problems, such as demand-side management in the smart grids \ncite{saad2012game}, charging/discharging of electric vehicles \ncite{grammatico2017dynamic} and formation control \ncite{lin2014distributed}, can be modelled as GNEPs.

The literature on NEPs and GNEPs mostly ignores the dynamics of individual agents, which may be problem in real-world multi-agent systems with non-negligible dynamics. The small portion of literature on NEPs and GNEPs with dynamical agents can be divided into two groups: algorithms that use the passivity properties of first-order methods and algorithms that use extremum seeking as a \textit{zero-order} method.

By using the passivity property in \ncite{gadjov2018passivity}, the authors design a control law that guarantees convergence to the Nash equilibrium of a multi-agent system with single-integrator dynamics over a network. In \ncite{romano2019dynamic}, the authors extend the result to the multi-integrator case. The network topology was considered to be non-time-varying in both cases. As a result, certain assumptions have to be met by the network graph. In \ncite{de2019distributed}, the authors extend the results of \ncite{gadjov2018passivity} by designing a network weight adaptation scheme. In \ncite{mattia2019dynamic}, a controller was proposed which guarantees convergence to a generalized Nash equilibrium of a multi-agent system with single and double integrator dynamics over a network. 

Most prominently, extremum seeking was used for Nash equilibrium seeking in \ncite{frihauf2011nash} where it was proven that the extremum seeking control, under certain conditions on the individual cost functions, will converge to the neighborhood of the Nash equilibrium for general nonlinear agents. In \ncite{liu2011stochastic}, it is proven that the use of stochastic perturbation signals also admits the convergence to the neighborhood of the Nash equilibrium. The authors in \ncite{poveda2017framework} propose a framework for synthesis of a hybrid controller which could also be used for NEPs with nonlinear agents. All of the mentioned extremum seeking controllers are based on  \ncite{krstic2000stability}.

\subsubsection{Contribution} Motivated by the recent research interest in NEPs, we adapt the extremum seeking  controller proposed in \ncite{guay2017proportional}, \ncite{guay2018distributed}. Specifically, our contributions are the following:
\begin{itemize}
    \item We extend a known proportional-integral extremum seeking control scheme to strongly monotone NEPs for a multi-agent linear systems and we prove a practical convergence to a Nash equilibrium.
    \item We numerically observe an improved performance with respect to \ncite{frihauf2011nash}, as smaller amplitudes and frequencies of the sinusoidal perturbations signals are needed for a comparable rate of convergence.
\end{itemize}

\subsubsection{Notation} $\mathbb{R}$ denotes the set of real numbers. For a matrix $A \in \mathbb{R}^{n \times m}$, $A^\top$ and $\|A\|$ denote its transpose and maximum singular value respectively. For vectors $x, y \in \mathbb{R}^{n}$, $x^\top y$ and $\|x \|$ denote the Euclidean inner product and norm, respectively. Given $N$ vectors $x_1, \dots, x_N$, possibly of different dimensions, $\bf{x} \coloneqq \left[ x_1^\top, \dots, x_N^\top \right]^\top $ and for each $i = 1, \dots, N$, $\bf{x}_{-i} \coloneqq \left[ x_1^\top, \dots,  x_{i -1}^\top,  x_{i + 1}^\top, \dots, x_N^\top \right]^\top $. $\operatorname{diag}\left(A_{1}, \ldots, A_{N}\right)$ denotes the block diagonal matrix with $A_{1}, \dots, A_{N}$ on its diagonal. For a function $v: \mathbb{R}^{n} \times \mathbb{R}^{m}  \rightarrow \mathbb{R}$ differentiable in the first argument, we denote the partial gradient vector as $\nabla_x v(x, y) \coloneqq \left[\frac{\partial v(x, y)}{\partial x_{1}}^\top, \ldots, \frac{\partial v(x, y)}{\partial x_{N}}^\top\right]^\top \in \mathbb{R}^{n}$. For a maping $v: \mathbb{R}^{n}  \rightarrow \mathbb{R}^m$, we denote the set of zeros as $\mathrm{zer}(V) \coloneqq \{ x \in \mathrm{dom}(A) | \bf{0}_{m} \in V(x) \}$.
\section{Problem setup}
We consider a multi-agent system with $N$ agents indexed by $\mathcal{I} = \{1, 2, \dots, N\}$, each with the following dynamics:
\begin{subequations}{}
\begin{align}
    \dot{x}_i &= -x_i+B_i u_i \label{sistem_i}\\
    y_i &= h_i(x_i, \bf{x}_{-i}), \label{cost_i}
\end{align}{}
\end{subequations}where $x_i \in \mathbb{R}^{n_i}$ is the state vector, $u_i \in \mathbb{R}^{m_i}$ is the control input,  $y_i \in \mathbb{R}$ is the output variable which evaluates the cost function $h_{i}: \mathbb{R}^{n_i} \times \mathbb{R}^{n_{-i}} \rightarrow \mathbb{R}$. Let us also define $n \coloneqq \sum n_i$, $n_{-i} \coloneqq \sum_{j \neq i} n_j$ and $m \coloneqq \sum m_i$.

\emph{{Standing Assumption 1 (Regularity)}}
\\
For each $i \in \mathcal{I}$, the function $h_i$ in \nref{cost_i} is differentiable in $x_i$ and its partial gradient $\nabla_{x_i} h_i$ is Lipschitz continuous in $x_i$ and $\bf{x}_{-i}$. \endassum

A common assumption amongst the extremum seeking literature (for example \ncite{krstic2000stability}, \ncite{guay2017proportional}, \ncite{poveda2017framework}) is the existence of the steady-state mapping, which tells us to which state(s) the system converges when a constant input is applied. For our subsystems \nref{sistem_i}, for each $i \in \mathcal{I}$, there exists a mapping 
\begin{align}
\pi(\bf{u}) \coloneqq \left[ \begin{array}{c}{\pi_1(u_1)} \\ \vdots \\ {\pi_N(u_N)}
\end{array}   \right]  =  \left[ \begin{array}{c}{B_1 u_1}  \\ \vdots \\{B_N u_N}
\end{array}   \right]\label{steady state mapping}
\end{align}
such that for every $i \in \mathcal{I}$,
\begin{align}
    \pi_i(u_i)  = B_i u_i.
\end{align}{}Let us also define
\begin{align}
\pi_{-i}(\bf{u}_{-i}) \coloneqq \left[ \begin{array}{c}{\pi_1(u_1)} \\ \vdots \\ {\pi_{i-1}(u_{i-1})} \\ {\pi_{i + 1}(u_{i + 1})} \\ \vdots \\{\pi_N(u_N)}
\end{array}   \right].
\end{align}

In this paper, we assume that the goal of each agent is to minimize its own steady-state cost function, i.e.,
\begin{align}
    \min _{u_{i} \in \mathbb{R}^{m_i}} h_{i}\left(\pi_i(u_i), \pi_{-i}(\bf{u}_{-i}) \right), \label{game}
\end{align}{}which depends on the inputs of some other agents as well. From a game-theoretic perspective, we consider the problem to compute a Nash equilibrium (NE). 

\begin{definition}[Nash equilibrium]
A collective input $\bf{u}^*$ is a NE of the game \nref{game} if for all $i \in \mathcal{I}$
\begin{equation*}
    \quad h_{i}\left(\pi_i(u_i^*), \pi_{-i}(\bf{u}_{-i}^*) \right) \leq \inf_{u_i \in \mathbb{R}^{n_i}} h_{i}\left(\pi_i(u_i), \pi_{-i}(\bf{u}_{-i}^*) \right). \qed\label{nash eq}
\end{equation*}{}
\end{definition}
In plain words, a set of inputs is a NE if no agent can improve its steady-state cost function by unilaterally changing its input. 

Since for all $i \in \mathcal{I}$, the steady-state cost functions are differentiable in $u_i$, it follows from Theorem 16.3 in \ncite{bauschke2011convex} that a collective vector $\bf{u}^*$ is a NE if and only if
\begin{equation}
    \nabla_{u_i} h_{i}\left(\pi_i(u_i^*), \pi_{-i}(\bf{u}_{-i}^*) \right) = 0. \label{gradijenti cijena}
\end{equation}{}
In view of \nref{gradijenti cijena}, we can stack all of the partial gradients into a single vector and form the so-called pseudo-gradient mapping of the steady-state cost functions:
\begin{align}
F(\bf{u}) \coloneqq \left[ \begin{array}{c}{\nabla_{u_1} h_1 \left( \pi_1(u_1), \pi_{-i}(\bf{u}_{-1}) \right)} \\ \vdots \\ {\nabla_{u_N} h_N \left( \pi_N(u_N), \pi_{-N}(\bf{u}_{-N}) \right)}
\end{array}   \right].  
\label{pseudogradient steady cost}
\end{align}Therefore, by \nref{gradijenti cijena} and \nref{pseudogradient steady cost}, we note that the problem of finding a Nash equilibrium of the game in \nref{game} is equivalent to finding $\bf{u}^*$ such that $F(\bf{u}^*) = 0$, which is the problem of finding a zero of $F$ in \nref{pseudogradient steady cost}, $\bf{u}^* \in \operatorname{zer}(F)$.

A relatively standard assumption in modern game theory literature \ncite{yu2017distributed}, \ncite{yi2019operator} is strong monotonicity of the pseudo-gradient:

\emph{{Standing Assumption 2 (Strong monotonicity)}}\\
The mapping $F$ in \nref{pseudogradient steady cost} is strongly monotone, i.e., 
\begin{align}
    (F(\bf{u}) - F(\bf{v}))^\top(\bf{u} - \bf{v}) \geq \mu \| \bf{u} - \bf{v} \|^2,
\end{align}{}for all $(\bf{u}, \bf{v}) \in \mathbb{R}^{2m}$, for some $\mu > 0$.

Let us also define the pseudo-gradient with respect to $x$ of the cost functions 
\begin{align}
    F_{\textup{x}}(\bf{x}) \coloneqq \left[ \begin{array}{c}{\nabla_{x_1} h_1 \left(x_1, \bf{x}_{-1} \right)} \\ \vdots \\ {\nabla_{x_N} h_N \left( x_N, \bf{x}_{-N} \right)}
\end{array}   \right]. \label{pseudogradient cost}
\end{align}{}

We note that, in general, monotonicity of the pseudo-gradient of the cost function $F_{\textup{x}}(\bf{x})$ does not imply monotonicity of the pseudo-gradient of the steady-state cost function $F(\bf{u})$. 

\section{Integral Nash equilibrium seeking control}
In this section, we propose two control schemes for Nash equilibrium seeking. The first one is designed for the full-information case. In other words, agents have perfect information of the actions of other agents and some additional information which will be described in more detail in the next subsection. The second is data-driven, i.e. designed for the case when the agents have access to the cost output. 

\subsection{Full-information case}
In the simplest case, we assume that every agent knows the analytic expression of its partial gradient and has access to the inputs of the other agents. The integral Nash equilibrium control in the next subsection will approximate the gradient and inputs of the other players and use the approximations in the same control law as the full-information case. Our proposed control law is inspired by the extremum seeking control in \ncite{guay2017proportional}, \ncite{guay2018distributed}:
\begin{align}
\forall i \in \mathcal{I}:\ \dot{u}_i = -\frac{1}{\tau_i}B_i^{\top} \nabla_{x_i} h_i(x_i, \bf{x}_{-i}) \label{controller full info}
\end{align}{}
or in collective vector form
\begin{align}
    \dot{\bf{u}} = -\bf{\tau}^{-1} \bf{B}^\top F_{\textup{x}}(\bf{x}),    \label{controller full info collective}
\end{align}{}where $\bf{B} \coloneqq \mathrm{diag}(B_1, \dots, B_N)$ and $\bf{\tau} \coloneqq \operatorname{diag}(\tau_{1}, \dots, \tau_{N})$. Unlike \ncite{guay2017proportional}, we do not use the proportional part, as it does not help with the convergence to the Nash equilibrium.
\begin{thm}[\textit{Convergence to Nash equilibrium}]\hfill
Consider the multi-agent system with dynamics \nref{sistem_i} and control law in \nref{controller full info}-\nref{controller full info collective} and let ($\bf{x}(t), \bf{u}(t)$) be its closed-loop solution. Let the Standing Assumptions hold. Then, there exists $\tau^*$ such that if $\min_{i \in \mathcal{I}} \tau_i \geq \tau^*$, then the pair $(\bf{x}(t), \bf{u}(t))$ converges to $(\bf{x}^*, \bf{u}^*) = (\pi(\bf{u}^*), \bf{u}^*)$, where $\bf{u}^*$ is a Nash equilibrium of the game in \nref{game}. \endassum
\end{thm}{}
\begin{pf}
See Appendix A.\hfill $\blacksquare$
\end{pf}{}

\subsection{Limited information case}
In the limited information case, we consider that the agents have access to the cost output. We emphasize that they neither know the actions of other agents, nor they know the analytic expressions of their partial gradients. This is the standard setup used in extremum seeking (\ncite{krstic2000stability}, \ncite{guay2017proportional}, \ncite{poveda2017framework} among others).

The extremum seeking control proposed by \ncite{guay2017proportional} assumes that the cost function of the system has a strong relative degree of value one. This means that the first derivative of the cost function has a direct influence on the input to the system. In the case of multi-agent systems, where the cost functions do not depend only on the states of their agent but also of the others, we make an analogous assumption:

\emph{{Assumption 1 (Degree of the output)}}\\
For every $i \in \mathcal{I}$, $\nabla_{x_i}h_i(x_i, \bf{x}_{-i})^\top B_i \neq 0 $ for all $(x_i, \bf{x}_{-i}) \in \mathbb{R}^n \setminus \{\bf{x}^*\}$.\endassum

Let us first evaluate the derivative of the cost functions:  
\begin{align}
\dot{y}_i &= -\sum_{j = 1}^N \nabla_{x_j} h_i(\bf{x})^\top + \sum_{j \neq i}^N \nabla_{x_j} h_i(\bf{x})^\top B_ju_j \nonumber \\ &+ \nabla_{x_i} h_i(\bf{x})^\top B_iu_i, \label{semiperfect y derivative}
\end{align}{}
and introduce the following variables:
\begin{align}
  \theta_{i}^0 &\coloneqq -\sum_{j = 1}^N \nabla_{x_j} h_i(\bf{x})^\top + \sum_{j \neq i}^N \nabla_{x_j} h_i(\bf{x})^\top B_j u_j,  \label{theta_0}\\
  \theta_{i}^1 &\coloneqq \nabla_{x_i} h_i(\bf{x})^\top B_i.\label{theta_1}
\end{align}{}The variable $\theta_{i}^{0}$ measures the effect of the autonomous dynamics of agent $i$ on its cost function and the effects of inputs of the other agents. The variable $\theta_{i}^{1}$ measures the effect of the input of agent $i$ on the cost output $y_i$. By substituting the previous variable definitions into \nref{semiperfect y derivative},  the derivative of the cost output reads as
\begin{align}
    \dot{y}_i = \theta_{i}^{0} + \theta_{i}^{1}u_i = [1, u_i^\top]\theta_i,\label{parametrizacija y}
\end{align}{}where $\theta_i = [\theta_{i}^{0}, \theta_{i}^{1\top}]$. Note that $\theta_{i}^{1}$ is proportional to the right-hand side in \nref{controller full info}. To estimate the local $\theta_{i}^{0}$ and $\theta_{i}^{1}$, we use a time-varying parameter estimation approach such as the one proposed in \ncite{guay}. Let us provide a basic intuition.\\
Let $\hat{y}_i$ and $\hat{\theta}_i$ be estimations of the output $y_i$ and the variable $\theta_i$ respectively and let $e_i = y_i - \hat{y}_i$ be the estimation error. Then, the estimator model of \nref{parametrizacija y} for agent $i$ is given by 
\begin{align}
    \dot{\hat{y}}_i&=[1, u_i^\top] \hat{\theta}_i+K_i e_i+c_i^{\top} \dot{\hat{\theta}}_i, \label{estim1}
\end{align}
where $K_i$ is a free design parameter. Note that the first two terms on the right-hand side resemble high gain observer schemes. As the structure of the problem does not allow the use of high gain observers, it is necessary to introduce some other dynamics into the estimation. This is the primary role of the third term in \nref{estim1}. Therefore, the dynamics of $c_i(t)$ are choosen as 
\begin{align}
\dot{c}_i^{\top}&=-K_i c_i^{\top}+[1, u_i^\top]. \label{ci dyn}
\end{align}
Let us also introduce an auxiliary variable $\eta_i$, with dynamics $\dot{\eta}_i=-K_i \eta_i-c_i^{\top} \dot{\theta}$, and its estimate $\hat{\eta}_i$, with dynamics
\begin{align}{}
\dot{\hat{\eta}}_i&=-K_i \hat{\eta_i}.
\end{align}
The original parameter estimation law in \ncite{adetola2008finite} was designed for constant parameters, therefore $\eta_i = \hat{\eta_i}$ and was fully known. For time-varying parameters, we still want to use $\eta_i$, but the additional term $-c_i^{\top} \dot{\theta}$ does not allow for its calculation, since the rate of change of the parameters is unknown. This is why the estimate $\hat{\eta}_i$ is used in the parameter estimation law. It is also necessary to define a symmetric, positive definite scaling matrix variable $\Sigma_i \in \mathbb{R}^{m_i + 1 \times m_i + 1}$ with dynamics given by
\begin{align}
&\dot{\Sigma}_i=c_i c_i^{\top}-k_{i}^\textup{T} \Sigma_i+\sigma_i I 
&\Sigma_i(0) = \alpha_{i}^1,\label{estim2}
\end{align}{}where $k_{i}^\top$, $\sigma_i$ and $\alpha_{1}^i$ are free design parameters. We note that, originally, in \ncite{adetola2008finite}, $\dot{\Sigma}_i=c_i c_i^{\top}$, but this proved to be inconvenient in practical implementations, as the elements of $\Sigma_i$ grow unbounded. Instead, as in \nref{estim2}, dynamics of $\Sigma_i$ behave as a first-order system. The third term is added so that the matrix is always invertible.

Equations \nref{estim1}-\nref{estim2} form the parameter update law presented in \ncite{adetola2008finite}:
\begin{align}
\dot{\hat{\theta}}_i=\operatorname{\Pi}_{\Theta_i}\left(\hat{\theta}_i, \Sigma_i^{-1}(c_i(e_i-\hat{\eta_i})-\sigma_i \hat{\theta}_i)\right), \label{estimlast}
\end{align}{}where $\operatorname{\Pi}_{\Theta_i}(\hat{\theta}, v)$ denotes the projection of the vector $v$ onto the tangent cone of the set $\Theta_i$ at $\hat{\theta}$, as defined by Equation 2.14 in \ncite{nagurney2012projected}. This implies that if the starting value $\hat{\theta}_i (0)$ is in $\Theta_i$, so will $\hat{\theta}_i(t) $ for all $t$. 

We are finally ready to propose an integral decentralized Nash equilibrium seeking control law of the form 
\begin{align}
\forall i \in \mathcal{I}:\left\{\begin{array}{c}{u_i=\hat{u}_i+d_i(t)}  \\ {\dot{\hat{u}}_i=-\frac{1}{\tau_i} \hat{\theta}_{i}^{1}\ \ \ \ }  \end{array}\right.\label{controll aprx}
\end{align}
together with Equations \nref{estim1}--\nref{estimlast}. In the collective vector form, Equation \nref{controll aprx} read as
\begin{align}
\left\{\begin{array}{c}{\bf{u}=\hat{\bf{u}}+\bf{d}(t)}  \\ {\dot{\hat{\bf{u}}}=-\bf{\tau}^{-1} \hat{\theta}^{1}}  \end{array}\right. \label{collective input}
\end{align}{}As in \ncite{guay2017proportional}, for the parameter estimation scheme to converge, a persistency of excitation (PE) assumption  for every agent is introduced.

\smallskip
\emph{{Assumption 2 (Persistence of excitation)}}\\
For every $i \in \mathcal{I}$, there exist constants $\alpha_{i}^2$ and $T_i$ such that 
\begin{align}
&\int_{t}^{t+T_i} c_i(\tau) c_i(\tau)^\top d \tau \geq \alpha_{i}^2 I,
&\forall t>0,
\end{align}{}where $c_i(\tau)$ is the solution to \nref{ci dyn}.\endassum

We conclude the section with the main theoretical result of the paper, namely, the convergence of the closed-loop dynamics to a Nash equilibrium of the game.
\begin{thm}[\textit{Convergence to Nash equilibrium}]\hfill
Consider a multi-agent system with dynamics \nref{sistem_i} and control law \nref{estim1} -- \nref{estimlast}, \nref{controll aprx}, for all $i \in \mathcal{I}$ and let ($\bf{x}(t), \bf{u}(t)$) be its closed-loop solution. Let the Standing Assumptions and Assumptions 1, 2 hold, let $\pi$ be the steady-state mapping in \nref{steady state mapping} and let $D$ be the largest amplitude of the perturbation signals $\{d_i(t)\}_{i \in \mathcal{I}}$. Then, there exist gains $(K_i, k_{i}^\textup{T}, \sigma_i)_{i \in \mathcal{I}}$ and $\tau^*$ such that if $\min_{i \in \mathcal{I}}\tau_i \geq \tau^*$, then the pair $(\bf{x}(t), \bf{u}(t))$ converges towards the $\mathcal{O}(D^2)$ neighborhood of some $(\bf{x}^*, \bf{u}^*) = (\pi(\bf{u}^*), \bf{u}^*)$, where $\bf{u}^*$ is a Nash equilibrium of the game in \nref{game}. \endassum
\end{thm}{}
\begin{pf}
See Appendix B. \hfill $\blacksquare$
\end{pf}{}

\section{Simulation example}
Consider a system with three agents $i \in \{ 1, 2, 3\}$, having dynamics:
\begin{align}
   \dot{x}_i = -x_i + u_i.
\end{align}{}
The cost functions of agents are given by
\begin{align}
    y_1 &= 1.5(x_1 - 1)^2 + 1.5x_1x_2 + x_1x_3 \nonumber \\
    y_2 &= -2x_2x_1 + 1.5(x_2 - 2)^2 + x_2x_3  \nonumber \\
    y_3 &= -2.5x_3x_1 -x_3x_2 + 1.5(x_3 - 3)^2.
\end{align}{}
\begin{figure}
    \centering
    \includegraphics[width = \linewidth]{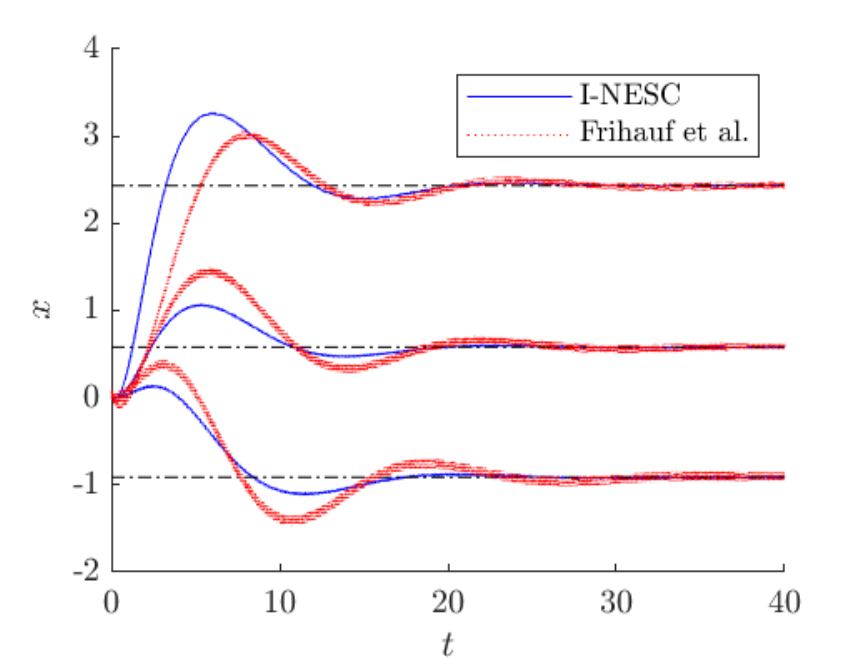}
    \caption{State trajectories of the three agents under I-NESC (solod blue) and \ncite{frihauf2011nash} (doted red)}
    \label{fig:comparison x}
\end{figure}{}
\begin{figure}
    \centering
    \includegraphics[width = \linewidth]{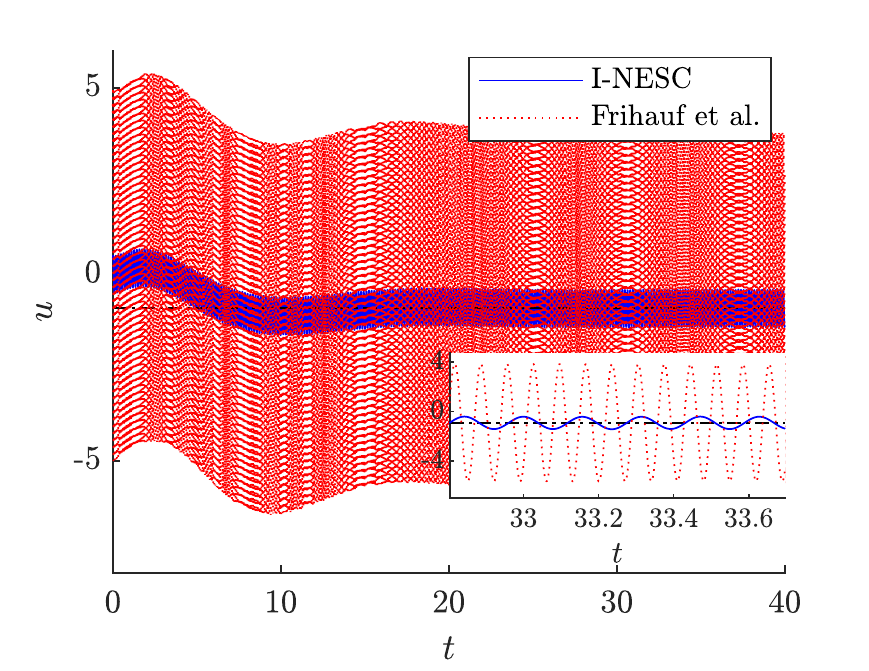}
    \caption{Input of the first agent under I-NESC (solid blue) and \ncite{frihauf2011nash} (doted red)} \label{fig:comparison_u}
\end{figure}{}Two types of controllers were simulated to have a comparison; the limited information controller proposed by this paper and the controller proposed in \ncite{frihauf2011nash} with additional low-pass and high-pass filters as in \ncite{krstic2000stability} in order to improve the performance. The latter can be described by the following equations
\begin{align} 
\dot{\eta}_i &=-\omega_{h}^i \eta_i+\omega_{h}^i y_i, \nonumber\\
\dot{\xi}_i &=-\omega_{l}^i \xi_i+\omega_{l}^i(y_i-\eta_i) A_i \sin (\omega_i t), \nonumber\\
\dot{\hat{u}}_i &=-k_i A_i \xi_i, \nonumber\\
u_i &=\hat{u}_i+A_i \sin (\omega_i t).
\end{align}
For our controller, the following parameters were chosen: $\sigma_1 = \sigma_2 = \sigma_3 = 10^{-6}$, $K_1 = K_2 = K_3 = 50$, $k_{T}^1 = k_{T}^2 = k_{T}^3 = 50$, $\alpha_I^1 = \alpha_I^2 = \alpha_I^3 = 0.1$, $\tau_{1} = 5$, $\tau_I^2 = 10$, $\tau_I^3 = 15$, $d_1(t) = \frac{1}{2}\sin{(40t)}$, $d_2(t) = \frac{1}{2}\sin{(50t)}$ and $d_3(t) = \frac{1}{2}\sin{(60t)}$. Initial states of $\bf{x}$, $\hat{\bf{u}}$, $\bf{c}$, $\bf{\theta}$ and $\bf{\eta}$ were set to zero. The parameters $K$, $k^\textup{T}$ and $\tau_I$ were initially chosen large enough to ensure stability. Then $\tau_I$ was decreased to speed up the convergence. Further decreases in $\tau_I$ were making the states oscillate; further decreasing of $K$ and $k^\textup{T}$ did not improve the performance of the algorithm.\\
For the Frihauf et al., the following parameters were chosen:  $\omega_h^1 = 180$, $\omega_h^2 = 200$, $\omega_h^3 = 220$, $\omega_l^1 = 45$, $\omega_l^2 = 50$, $\omega_l^3 = 55$, $\omega_1 = 90$, $\omega_2 = 100$, $\omega_3 = 110$, $k_1 = k_2 = k_3 = 0.5$ and $A_1 = A_2 = A_3 = 5$. The parameters were experimentally chosen in such a way that a fast convergence rate is obtained without increasing the amplitude of the sinusoidal perturbations to unreasonable values. 

The results of the numerical simulations can be seen in Figures \ref{fig:comparison x} and \ref{fig:comparison_u}. While the convergence speed of both algorithms is similar, the frequency and amplitude of the sinusoidal perturbation signals are much lower with our I-NESC law.

We remark that the proportional-integral extremum seeking control in \ncite{guay2017proportional} is capable of achieving a faster convergence speed than the algorithm in \ncite{krstic2000stability}. Our I-NESC algorithm has slower convergence compared to \ncite{guay2017proportional} because of the lack of the proportional part in the control law. Although for optimization problems it can greatly improve convergence speed, we observed that for NEPs even small proportional gain may cause instability. 
\section{Conclusion}
Nash equilibrium problems can be solved by an extremum seeking algorithm if the agents belong to a certain class of linear dynamics with strongly monotone and Lipschitz continuous game mapping. Extension to general linear or nonlinear systems and the extension to constrained states and inputs, is left for future work.

                                                   







\appendix
\section{Proof of Theorem 1}    
To prove the theorem, it is necessary to show that our subsystems, which under constant inputs stabilize to some equilibrium points, will not become unstable when a time-varying input is applied. Additionally, it is necessary that our collective input $\bf{u}$ converges to the Nash equilibrium $\bf{u}^*$. First, we will describe the behaviour of the subsystems to time-varying inputs.

Stability of equilibrium points $\pi_i(u_i)$ for every agent $i \in \mathcal{I}$ can be characterised by the following Lyapunov function:
\begin{align}
    V_i(x_i, u_i) &= \frac{\beta}{2}(x_i - \pi_i(u_i))^\top(x_i - \pi_i(u_i)), \label{lyapunov agent i} \nonumber \\
     &= \frac{\beta}{2}(x_i - B_i u_i)^\top (x_i - B_i u_i),
\end{align}{}
where $\beta > 0$. The derivative of \nref{lyapunov agent i} is equal to 
\begin{align}
    \dot{V}_i(x_i, u_i) &= \nabla_x V_i(x_i, u_i)^\top \dot{x}_i + \nabla_u V_i(x_i, u_i)^\top \dot{u}_i \nonumber \\
    &  = -\beta \|x_i - \pi_i ( u_i) \|^2 \nonumber \\
&- \beta (x_i - \pi_i ( u_i))\top B \dot{u}_i. \label{semiperfect izvod lyapunov i}
\end{align}{}Note that if $u_i$ is constant ($\dot{u}_i = 0$), equilibrium point $\pi_i(u_i)$ is uniformly globally exponentially stable (UGES). For very slow changes of the input $u_i$, we expect the subsystems will converge to a small neighborhood of equilibrium points $\pi_i(u_i)$, while for fast changes of the input $u_i$, we do not. We consider the controller in \nref{controller full info collective}. Its goal is to both estimate the Nash equilibrium input $\bf{u}^*$ and to preserve the stability of the subsystems. Therefore, we construct the following Lyapunov function candidate:
\begin{align}
    W(\bf{x}, \bf{u}) \coloneqq T(\bf{u}) + V(\bf{x} , \bf{u}) = \frac{1}{2 \tau_{\textup{min}}}\Tilde{\bf{u}}^{\top} \bf{\tau}\Tilde{\bf{u}} + \sum_{i = 1}^N V_i(x_i , u_i), \label{new_lypunov_i}
\end{align}{}where 
$\Tilde{\bf{u}} = \bf{u} - \bf{u}^*$ and $\tau_{\textup{min}} = \min\{\tau_{1}, \dots, \tau_{N}\}$. 

Now, we bound the derivative of $T$. By adding and subtracting $F(\bf{u})$ to \nref{controller full info collective}, $\dot{\bf{u}}$ reads as
\begin{align}
    \dot{\bf{u}} = -\bf{\tau}^{-1}F(\bf{u}) -\bf{\tau}^{-1} (\bf{B}^\top F_{\textup{x}}(\bf{x}) - F(\bf{u})). \label{contro law expanded}
\end{align}{}From \nref{new_lypunov_i} and \nref{contro law expanded}, we have
\begin{align}
    \dot{T}(\bf{x}, \bf{u}) = -\frac{1}{\tau_{\textup{min}}}\Tilde{\bf{u}}^\top F(\bf{u}) -\frac{1}{\tau_{\textup{min}}}\Tilde{\bf{u}}^{\top} (\bf{B} F_{\textup{x}}(\bf{x}) - F(\bf{u})).
\label{lypunov input}
\end{align}{}Since $F(\bf{u})$ is strongly monotone, it holds that
\begin{align}
    (F(\bf{u}) - F(\bf{u}^*))^{\top} (\bf{u} - \bf{u}^*) \geq \mu \|\bf{u} - \bf{u}^*\|^2. \label{strong mono}
\end{align}As $\bf{u}^* \in \mathrm{zer}(F)$, we have $F(\bf{u}^*) = 0$, thus \nref{strong mono} reads as
\begin{align}
    F(\bf{u})^{\top} (\bf{u} - \bf{u}^*)  = F(\bf{u})^{\top} \Tilde{\bf{u}} \geq \mu \|\bf{u} - \bf{u}^*\|^2. \label{strong mono 2}
\end{align}{}To bound the second term in \nref{lypunov input}, we use the following identity:
\begin{align}
    &\nabla_{u_i} h_i (\pi_i(u_i), \pi_{-i}(\bf{u}_{-i}))^\top \nonumber \\
    & = \nabla_{x_i} h_i (\pi_i(u_i), \pi_{-i}(\bf{u}_{-i}))^\top \nabla_{u_i}\pi_{i}(u_i) \nonumber \\
    & = \nabla_{x_i} h_i (\pi_i(u_i), \pi_{-i}(\bf{u}_{-i}))^\top B_i. \label{gradient relation}
\end{align}{}By using the relations \nref{pseudogradient steady cost}, \nref{pseudogradient cost}  and \nref{gradient relation}, it follows that:
\begin{equation}
    \bf{B}^\top F_{\textup{x}}(\pi(\bf{u})) = F(\bf{u}). \label{cost connection}
\end{equation}{}By exploiting \nref{strong mono 2} and \nref{cost connection}, from \nref{lypunov input} we have
\begin{align}
&\dot{T} \leq -\frac{\mu}{\tau_{\textup{min}}} \|\bf{u} - \bf{u}^*\|^2 \nonumber \\
&-\frac{1}{\tau_{\textup{min}}}\Tilde{\bf{u}}^{\top} (\bf{B}^\top F_{\textup{x}}(\bf{x}) - \bf{B}^\top F_{\textup{x}}(\pi(\bf{u}))).
\label{lypunov input 2}
\end{align}{}Since all of the functions are Lipschitz continuous, the right-hand side in \nref{lypunov input 2} can be upper bounded as follows:
\begin{align}
    \dot{T} \leq -\frac{\mu}{\tau_{\textup{min}}} \|\Tilde{\bf{u}}\|^2 +\frac{L}{\tau_{\textup{min}}}\|\Tilde{\bf{u}}\|\|\bf{x} - \pi(\bf{u}) \|, \label{bound T}
\end{align}
where $L > 0$ is the Lipschitz constant of the mapping $\bf{B}^\top \circ F_{\textup{x}}$. 

Now, we turn our attention to the full Lyapunov function candidate $W$. The derivative of $W$ can be now bounded as 
\begin{align}
    \dot{W}(\bf{x}, \bf{u}) &\leq  -\beta \|\bf{x}-\pi(\bf{u})\|^{2}  - \frac{\mu }{\tau_{\textup{min}}}\|\Tilde{\bf{u}}\|^2 \nonumber \\ & +\frac{L}{\tau_{\textup{min}}}\|\Tilde{\bf{u}}\|\|\bf{x} - \pi(\bf{u}) \| - \beta (\bf{x}- \pi(\bf{u}))\top \bf{B} \dot{\bf{u}}. \label{semi derivative W}
\end{align}
\\
To complete the proof, we bound the derivative of $V$ caused by the change of inputs:
\begin{align}
    -\beta (\bf{x}- \pi(\bf{u}))\top \bf{B} \dot{\bf{u}} \leq \beta \|\bf{B}\|\|\dot{\bf{u}}\|\|\bf{x}- \pi(\bf{u})\|. \label{} \label{bound derivate lyap}
\end{align}{}By using \nref{contro law expanded}, the norm of the derivative is bounded:
\begin{align}
    \|\dot{\bf{u}}\| 
    &\leq \frac{1}{\tau_{\textup{min}}} \|F(\bf{u}) \| + \frac{1}{\tau_{\textup{min}}} \|(\bf{B}^\top F_{\textup{x}}(\bf{x}) - F(\bf{u}))\|.
\end{align}{}
\\
Again, since all of the functions are Lipschitz continuous, the right-hand side of the previous equation can be bounded as follows
\begin{align}
    \|\dot{\bf{u}}\| \leq \frac{L_F}{\tau_{\textup{min}}} \|\tilde{\bf{u}} \| + \frac{L}{\tau_{\textup{min}}} \|\bf{x} - \pi(\bf{u})\|, \label{bound derivative}
\end{align}{}where $L_F > 0$ is the Lipschitz constant of $F$. 
By using the bounds \nref{bound derivate lyap} and\nref{bound derivative}, $\dot{W}$ can be bounded as follows:
\begin{align}
    \dot{W}(\bf{x}, \bf{u}) 
    &\leq -\left(\beta -\frac{L \beta \|\bf{B}\|}{\tau_{\textup{min}}}\right) \|\bf{x}-\pi(\bf{u})\|^{2} \nonumber \\ 
    & + \frac{L + L_F \beta \|\bf{B}\|}{\tau_{\textup{min}}} \|\tilde{\bf{u}} \|\| \bf{x} - \pi(\bf{u})\| - \frac{\mu }{\tau_{\textup{min}}}\|\Tilde{\bf{u}}\|^2 \nonumber\\
    & \leq  -\left[\begin{array}{c}{\left\|\Tilde{\bf{u}}\right\|} \\ {\left\|\bf{x}-\pi(\bf{u})\right\|}\end{array}\right]^\top M\left[\begin{array}{c}{\left\|\Tilde{\bf{u}}\right\|} \\ {\left\|\bf{x}-\pi(\bf{u})\right\|}\end{array}\right]
    \label{perfect derivative W} 
\end{align}where
\begin{align}
    M = \left[\begin{array}{cc}{\beta -\frac{L \beta \|\bf{B}\|}{\tau_{\textup{min}}}} & {-\frac{L + \beta \|\bf{B}\| L_F }{2\tau_{\textup{min}}}} \\ {-\frac{L + \beta \|\bf{B}\| L_F}{2\tau_{\textup{min}}}} & {\frac{\mu }{\tau_{\textup{min}}}}\end{array}\right]. \label{matica M}
\end{align}{}Thus, if
\begin{align}
    \tau_{\textup{min}} \geq \frac{(L  + \beta \|\bf{B}\| L_F)^2  + 4 L \beta \|\bf{B}\|}{4\beta \mu},
\end{align}{}then the matrix $M$ in \nref{matica M} is positive definite, which in turn implies that the Lyapunov derivative $\dot{W}$ is negative definite, which concludes the Lyapunov argument and in turn the proof. \hfill $\blacksquare$

\section{Proof of Theorem 2}         
The proof is similar to the full-information case proof, but unlike the full-information case, our inputs use the estimation of the $\theta_i$ variables, which affects stability. Let us consider a Lyapunov function candidate of the form 
\begin{align}
    L = W + V + T,
\end{align}{}where 
\begin{align}
    W(\tilde{\bf{\eta}}, \tilde{\bf{\theta}}) &=\sum_{i=1}^{N}\left(\frac{1}{2} \widetilde{\eta}_{i}^\top \widetilde{\eta}_{i}+\frac{1}{2} \widetilde{\theta}_{i}^\top \Sigma_{i} \widetilde{\theta}_{i}\right),\\
    V(\bf{x}, \hat{\bf{u}}) &= \sum_{i = 1}^N V_i(x_i, \hat{u}_i),\\
    T(\hat{\bf{u}}) &= \frac{1}{2\tau_{\textup{min}}}(\hat{\bf{u}}  - \bf{u}^*)^\top \bf{\tau}({\bf{u}}  - \bf{u}^*) = \frac{1}{2\tau_{\textup{min}}}\tilde{\bf{u}}^\top\bf{\tau}\tilde{\bf{u}}. \label{lyapnov T}
\end{align}{}Therefore, our Lyapunov function candidate consists of:
\begin{itemize}
    \item[(W)] Parameter estimation scheme\\
    \item[(V)] Local state-input Lyapnuov function\\
    \item[(T)] Nash equilibrium estimation error
\end{itemize}

\subsubsection{Parameter estimation scheme} \hfill

We bound the time derivative of the W function similarly to \ncite{guay} and \ncite{guay2018distributed} with the only difference that we let each agent choose their own parameters ($\sigma_i, K_i, k_{i}^\textup{T}$). The Lyapunov derivative reads as follows:
\begin{align}
    \dot{W}(\tilde{\bf{\eta}}, \tilde{\bf{\theta}}) \leq& \sum_{i=1}^{p}\Bigg (-\widetilde{\eta}_{i}^\top\left(K_i-\frac{1}{2}-\frac{k_{1} \zeta_i}{2}\right) \widetilde{\eta}_{i} \nonumber \\ 
    & +\frac{1}{2 k_{1}} \dot{\theta}_{i}^\top \dot{\theta}_{i}+\frac{\gamma_{2i}}{2 k_{2}} \dot{\theta}_{i}^\top \dot{\theta}_{i}\frac{k_{i}^{T\prime} \gamma_{1i}}{2} \widetilde{\theta}_{i}^{\tau} \widetilde{\theta}_{i}+\frac{\sigma_i}{2} \theta_{i}^\top \theta_{i} \Bigg )\nonumber \\
    \leq& - k_\textup{a} \|\tilde{\bf{\eta}} \|^2 - k_\textup{b}\| \tilde{\bf{\theta}}\|^2 + k_\textup{c} \| \dot{\bf{\theta}} \|^2 + \frac{\sigma}{2}\|\bf{\theta}\|^2, \label{bound W}
\end{align}where $k_\textup{a} \coloneqq \min_i \left(K_i-\frac{1}{2}-\frac{k_{1} \zeta_i}{2}\right)$, $k_\textup{b} \coloneqq \min_i \left(\frac{k_{i}^{T\prime} \gamma_{1i}}{2}\right)$, $k_\textup{c} \coloneqq \max_i \left(\frac{1}{2k_1} + \frac{\gamma_{2i}}{2k_2}\right)$ and $\sigma \coloneqq\max_i\sigma_i$.
\subsubsection{Local state-input Lyapnuov function}\hfill\\
As stated in the previous section (\nref{lyapunov agent i} and \nref{semiperfect izvod lyapunov i}), equilibrium points $\pi_i(u_i)$ are UGES when a constant input is applied. 
Stability of equilibrium points $\pi_i(\hat{u}_i)$ for time-varying inputs, for every agent $i \in \mathcal{I}$ can be characterised by the following Lyapunov function:
\begin{align}
    V_i(x_i, u_i) = \frac{\beta}{2}(x_i - \pi_i(\hat{u}_i))^\top (x_i - \pi_i(\hat{u}_i))^\top \label{lyapunov agent i estimate}  &, \beta > 0.
\end{align}{}The derivative of the Lyapunov function $V_i(x_i,u_i)$ is
\begin{align}
\dot{V}_i(x_i, u_i) =& -\beta \|x_i - \pi_i(\hat{u}_i)\|^2  + \beta (x_i - \pi_i(\hat{u}_i))^\top B_i d_i(t) \nonumber \\
&- \beta (x_i - \pi_i(\hat{u}_i))\top B \dot{\hat{u}}_i. \label{izvod lyapunov i estimate}
\end{align}{}Note that the derivative has three addends; the first one is equal to the complete derivative of the Lyapunov function in the case of constant inputs, the amplitude of the second component is proportional to the amplitude of the perturbations and the amplitude of the third component is equal to the amplitude of the derivative of the input $u_i$.

To bound the third component, we need to bound $\dot{\hat{u}}_i$, hence $\dot{\hat{\bf{u}}}$, which reads as
\begin{align}
    \dot{\hat{\bf{u}}}&=-\bf{\tau}^{-1} ({\bf{\theta}}^{1} + \tilde{\bf{\theta}}^{1} ) \nonumber \\
    &=-\bf{\tau}^{-1} (\bf{B}\top F_{\textup{x}}(\bf{x}) + \tilde{\bf{\theta}}^{1} ).
\end{align}
By using the same argument as in \nref{bound derivative}, it follows that
\begin{align}
    \dot{\hat{\bf{u}}} &\leq \frac{L_F}{\tau_{\textup{min}}} \|\tilde{\bf{u}} \| + \frac{L}{\tau_{\textup{min}}} \|\bf{x} - \pi(\hat{\bf{u}})\| +\frac{1}{\tau_{\textup{min}}} \|\tilde{\bf{\theta}}^{1} \| . \label{bound u changing}
\end{align}
By using the previous equation, it is possible to bound the second addend in \nref{izvod lyapunov i estimate}: 
\begin{align}
    &\beta (\bf{x} - \pi(\bf{u}))\top \bf{B} \dot{\hat{\bf{u}}} \leq \frac{L_F \beta \|\bf{B}\|}{\tau_{\textup{min}}} \|\tilde{\bf{u}} \|\|\bf{x} - \pi(\hat{\bf{u}})\|\nonumber \\ & + \frac{L \beta \|\bf{B}\|}{\tau_{\textup{min}}} \|\bf{x} - \pi(\hat{\bf{u}})\|^2 +\frac{\beta \|\bf{B}\|}{\tau_{\textup{min}}} \|\tilde{\bf{\theta}}^{1} \|\|\bf{x} - \pi(\hat{\bf{u}})\|.
\end{align}{}
Therefore, the derivative of $V$ can be bounded as
\begin{align}
      \dot{V}(\bf{x}, {\bf{u}}) &\leq - \beta\|\bf{x}-\pi(\hat{\bf{u}})\|^{2} + \frac{L_F \beta \|\bf{B}\|}{\tau_{\textup{min}}} \|\tilde{\bf{u}} \|\|\bf{x} - \pi(\hat{\bf{u}})\| \nonumber \\ &+ \frac{L \beta \|\bf{B}\|}{\tau_{\textup{min}}} \|\bf{x} - \pi(\hat{\bf{u}})\|^2 +\frac{\beta \|\bf{B}\|}{\tau_{\textup{min}}} \|\tilde{\bf{\theta}}^{1} \|\|\bf{x} - \pi(\hat{\bf{u}})\| \nonumber \\ & -\beta (\bf{x} - \pi(\hat{\bf{u}}))\top \bf{B}\bf{d}(t). \label{semiimperfect bound V}
\end{align}{}The last term can be bounded by using Cauchy–-Bunyakovsky–-Schwarz inequality to obtain
\begin{align}
     \dot{V}(\bf{x}, {\bf{u}})  \leq& - \beta\|\bf{x}-\pi(\hat{\bf{u}})\|^{2} + \frac{L_F \beta \|\bf{B}\|}{\tau_{\textup{min}}} \|\tilde{\bf{u}} \|\|\bf{x} - \pi(\hat{\bf{u}})\| \nonumber \\ &+ \frac{L \beta \|\bf{B}\|}{\tau_{\textup{min}}} \|\bf{x} - \pi(\hat{\bf{u}})\|^2 +\frac{\beta \|\bf{B}\|}{\tau_{\textup{min}}} \|\tilde{\bf{\theta}}^{1} \|\|\bf{x} - \pi(\hat{\bf{u}})\| \nonumber \\ &  +\beta \|\bf{B}\| \|\bf{x} - \pi(\hat{\bf{u}})\|\|\bf{d}(t)\|. \label{semiperfect bound V}
\end{align}{}Finally, the last two non-quadratic terms can be bounded with the inequality $\|x\|\|y\| \leq \frac{1}{2k}\|x\|^2 + \frac{k}{2}\|y\|^2$ to finally conclude the desired bound:
\begin{align}
    \dot{V}\leq& - \left(\beta  - \frac{L \beta \|\bf{B}\|}{\tau_{\textup{min}}} - \frac{\beta \|\bf{B}\|}{2\tau_{\textup{min}}k_3} - \frac{\beta \|\bf{B}\|}{2k_4}\right)\|\bf{x}-\pi(\hat{\bf{u}})\|^{2} \nonumber \\ 
    &+ \frac{L_F \beta \|\bf{B}\|}{\tau_{\textup{min}}} \|\tilde{\bf{u}} \|\|\bf{x} - \pi(\hat{\bf{u}})\| +  \frac{\beta \|\bf{B}\| k_3}{2\tau_{\textup{min}}} \|\tilde{\bf{\theta}}^{1} \|^2 \nonumber \\ 
    &+ \frac{\beta \|\bf{B}\| k_4}{2} \|\bf{d}(t)\|^2.    \label{perfect bound V}\end{align}{}\subsubsection{Nash equilibrium estimation error}                      \hfill\\
The use of parameter estimation also has an influence on the Nash equilibrium estimation error. The derivative of \nref{lyapnov T} is equal to
\begin{align}
    \dot{T}(\bf{x}, \hat{\bf{u}}) = -\bf{\tau}^{-1}\tilde{\bf{u}} (\bf{B}\top F_{\textup{x}}(\bf{x}) + \tilde{\bf{\theta}}_{i}^{1} ).
\end{align}{}
By the same method as in \nref{bound T}, \nref{perfect bound V}, it follows that
\begin{align}
    \dot{T}(\hat{\bf{u}}) \leq& -\frac{\mu}{\tau_{\textup{min}}} \|\Tilde{\bf{u}}\|^2 +\frac{L}{\tau_{\textup{min}}}\|\Tilde{\bf{u}}\|\|\bf{x} - \pi(\hat{\bf{u}}) \| + \frac{1}{\tau_{\textup{min}}}\|\tilde{\bf{u}}\|\|\tilde{\bf{\theta}}^{1}\| \nonumber \\
    \leq& -\frac{\mu}{\tau_{\textup{min}}} \|\Tilde{\bf{u}}\|^2 +\frac{L}{\tau_{\textup{min}}}\|\Tilde{\bf{u}}\|\|\bf{x} - \pi(\hat{\bf{u}}) \| + \frac{1}{2\tau_{\textup{min}}k_5}\|\tilde{\bf{u}}\|^2 \nonumber \\
    & + \frac{k_5}{2\tau_{\textup{min}}}\|\tilde{\bf{\theta}}^{1}\|^2. \label{bound T estimation}
\end{align}{}\subsubsection{The full Lyapunov candidate} \hfill \\
Now, with the bounds \nref{bound W}, \nref{semiperfect bound V} and \nref{bound T estimation}, the time derivative of the full Lyapunov candidate function is bounded as follows:
\begin{align}
    \dot{L} \leq& - k_\textup{a} \|\tilde{\bf{\eta}} \|^2 - k_\textup{b}\| \tilde{\bf{\theta}}^{0}\|^2 - \left(\frac{\mu}{\tau_{\textup{min}}} -  \frac{1}{2\tau_{\textup{min}}k_5}\right)\|\Tilde{\bf{u}}\|^2 \nonumber \\ &- \left(\beta  - \frac{L \beta \|\bf{B}\|}{\tau_{\textup{min}}}- \frac{\beta \|\bf{B}\|}{2\tau_{\textup{min}}k_3}  - \frac{\beta \|\bf{B}\|}{2k_4}\right)\|\bf{x}-\pi(\hat{\bf{u}})\|^{2}  \nonumber \\ & - \left(k_\textup{b} - \frac{k_5}{2\tau_{\textup{min}}} - \frac{\beta \|\bf{B}\| k_3}{2\tau_{\textup{min}}}\right)\| \tilde{\bf{\theta}}^{1}\|^2 + k_\textup{c} \| \dot{\bf{\theta}} \|^2 + \frac{\sigma}{2}\|\bf{\theta}\|^2  \nonumber \\ &  + \frac{L_F \beta \|\bf{B}\| + L}{\tau_{\textup{min}}} \|\tilde{\bf{u}} \|\|\bf{x} - \pi(\hat{\bf{u}})\|+ \frac{\beta \|\bf{B}\| k_4}{2} \|\bf{d}(t)\|^2. \label{semiperfect L bound}
\end{align}{}We are left with determining bounds on $\|\bf{\theta}\|$ and $\|\dot{\bf{\theta}}\|$. Since all of the considered functions (and their composition) in \nref{sistem_i}, \nref{theta_0}, \nref{theta_1} and \nref{collective input} are Lipschitz continuous, it follows that
\begin{align}
    \|\bf{\theta}\|^2 &\leq L_1 \|\bf{x} - \pi(\hat{\bf{u}})\|^2 + L_2 \|\tilde{\bf{u}}\|^2 \label{bound theta}\\
    \|\dot{\bf{\theta}}\|^2 &\leq L_3 \|\bf{x} - \pi(\hat{\bf{u}})\|^2 + L_4 \|\tilde{\bf{u}}\|^2 \label{bound theta dot},
\end{align}{}for some $L_1, L_2, L_3, L_4 > 0$. Substituting \nref{bound theta} and \nref{bound theta dot}  into \nref{semiperfect L bound}, we obtain
\begin{align}
        \dot{L} \leq& - \left(\frac{\mu}{\tau_{\textup{min}}} -  \frac{1}{2\tau_{\textup{min}}k_5} - \frac{L_2 \sigma}{2} - k_\textup{c} L_4\right)\|\Tilde{\bf{u}}\|^2 - \bigg (\beta  - \frac{L_1 \sigma}{2}\nonumber \\ 
        &- \frac{L \beta \|\bf{B}\|}{\tau_{\textup{min}}}- \frac{\beta \|\bf{B}\|}{2\tau_{\textup{min}}k_3}  - \frac{\beta \|\bf{B}\|}{2k_4} -  k_\textup{c} L_3  \bigg )\|\bf{x}-\pi(\hat{\bf{u}})\|^{2}  \nonumber \\ 
        & - \left(k_\textup{b} - \frac{k_5}{2\tau_{\textup{min}}} - \frac{\beta \|\bf{B}\| k_3}{2\tau_{\textup{min}}}\right)\| \tilde{\bf{\theta}}^{1}\|^2   - k_\textup{a} \|\tilde{\bf{\eta}} \|^2 - k_\textup{b}\| \tilde{\bf{\theta}}^{0}\|^2  \nonumber \\ 
        &+ \frac{L_F \beta \|\bf{B}\| + L}{\tau_{\textup{min}}} \|\tilde{\bf{u}} \|\|\bf{x} - \pi(\hat{\bf{u}})\| + \frac{\beta \|\bf{B}\| k_4}{2} \|\bf{d}(t)\|^2. \label{perfect L bound}
\end{align}{}The task at hand now is to prove that there exist parameters $K$, $k^\textup{T}$ and $\tau_{\textup{min}}$ such that the right-hand side in \nref{perfect L bound}, apart from the term with $\|\bf{d}(t)\|^2$, is always negative definite. The proof goes by the same lines as in \ncite{guay}, \ncite{guay2018distributed}. \\
Consider the following reformulation of \nref{perfect L bound}:
\begin{align}
        \dot{L} \leq& - k_\textup{a} \|\tilde{\bf{\eta}} \|^2 - k_\textup{b}\| \tilde{\bf{\theta}}^{0}\|^2 - \left(- \frac{L_2 \sigma}{2} - k_\textup{c} L_4\right)\|\Tilde{\bf{u}}\|^2 \nonumber \\ 
        &- \left(- \frac{\beta \|\bf{B}\|}{2k_4} - \frac{L_1 \sigma}{2} - k_\textup{c} L_3\right)\|\bf{x}-\pi(\hat{\bf{u}})\|^{2}  \nonumber \\ 
        & -\left(k_\textup{b} - \frac{k_5}{2\tau_{\textup{min}}} - \frac{\beta \|\bf{B}\| k_3}{2\tau_{\textup{min}}}\right)\| \tilde{\bf{\theta}}^{1}\|^2  + \frac{\beta \|\bf{B}\| k_4}{2} \|\bf{d}(t)\|^2 \nonumber \\
        & -\left[\begin{array}{c}{\left\|\Tilde{\bf{u}}\right\|} \\ {\left\|\bf{x}-\pi(\bf{u})\right\|}\end{array}\right]^\top \bf{M}\left[\begin{array}{c}{\left\|\Tilde{\bf{u}}\right\|} \\ {\left\|\bf{x}-\pi(\bf{u})\right\|}\end{array}\right]\label{perfect L bound with matrix},
\end{align}{}where 
\begin{align}
    \bf{M} = \left[\begin{array}{cc}{\beta -\frac{L \beta \|\bf{B}\|}{\tau_{\textup{min}}} - \frac{\beta \|\bf{B}\|}{2\tau_{\textup{min}} k_3}} & {-\frac{L + \beta \|\bf{B}\| L_F }{2\tau_{\textup{min}}}} \\ {-\frac{L + \beta \|\bf{B}\| L_F}{2\tau_{\textup{min}}}} & {\frac{2\mu - 1/k_5 }{2\tau_{\textup{min}}}}\end{array}\right].
\end{align}{}The parameter $k_3$ can be chosen arbitrarily, while $k_5$ has to be chosen such that lower diagonal element in $\bf{M}$ is positive, i.e. $2\mu - 1/k_5 > 0$. Also, in order for $\bf{M}$ to be positive definite, the following condition must be satisfied
\begin{align}
    &v_1(L, \beta \|\bf{B}\|, L_F, k_3, k_5, \beta, \mu)  = \nonumber \\
   &\frac{(L  + \beta \|\bf{B}\| L_F)^2  + 2 L \beta \|\bf{B}\| + \beta \|\bf{B}\|/k_3}{2\beta (2\mu - 1/k_5)} \leq  \tau_{\textup{min}}.
\end{align}{}From Equation \nref{perfect L bound with matrix}, it can be also concluded that 
\begin{align}
    \tau_{\textup{min}} \geq \frac{Lk_3 + k_5}{k_\textup{b}} = v_2(L, k_3, k_5, k_\textup{b})
\end{align}{}Therefore, $\tau_{\textup{min}} \geq \max(v_1, v_2)$. Let $\lambda$ be the smallest eigenvalue of matrix $\bf{M}$. The inequality \nref{perfect L bound with matrix} can be reformulated as
\begin{align}
        \dot{L} \leq& - k_\textup{a} \|\tilde{\bf{\eta}} \|^2 - k_\textup{b}\| \tilde{\bf{\theta}}^{0}\|^2 - \left(\lambda- \frac{L_2 \sigma}{2} - k_\textup{c} L_4\right)\|\Tilde{\bf{u}}\|^2 \nonumber \\ &- \left(\lambda- \frac{\beta \|\bf{B}\|}{2k_4} - \frac{L_1 \sigma}{2} - k_\textup{c} L_3\right)\|\bf{x}-\pi(\hat{\bf{u}})\|^{2}  \nonumber \\ & - \left(k_\textup{b} - \frac{k_5}{2\tau_{\textup{min}}} - \frac{\beta \|\bf{B}\| k_3}{2\tau_{\textup{min}}}\right)\| \tilde{\bf{\theta}}^{1}\|^2  + \frac{\beta \|\bf{B}\| k_4}{2} \|\bf{d}(t)\|^2. \nonumber \\
\end{align}{}The parameters $\sigma$, $k_\textup{c}$ must be chosen small enough, while $k_4$ must be large enough such that the following equations hold true: 
\begin{align}
    0 &< \lambda- \frac{L_2 \sigma}{2} - k_\textup{c} L_4, \nonumber \\
    0 &< \lambda- \frac{\beta \|\bf{B}\|}{2k_4} - \frac{L_1 \sigma}{2} - k_\textup{c} L_3.
\end{align}{}The parameter $\sigma$ is a free design parameter; the parameter $k_\textup{c}$ can be made arbitrarily small by increasing the gains $K$ and $k^\textup{T}$ (or to be more precise $K_i$ and $k_{i}^\textup{T}$, see \ncite{guay} for more details). The parameter $k_4$ can be arbitrarily chosen. Therefore, it is possible to choose the controller parameters $\sigma$, $k^\textup{T}$ and $K$ such that all of the constants that multiply the squares of the norms in \nref{perfect L bound with matrix} (except for $\bf{d}(t)$) are positive. \\
Next, we consider the Lyapunov functions of the subsystems in \nref{lyapunov agent i estimate}, the bounds on matrices $\Sigma_i$ and the quadratic elements of the Lyapunov function candidate $L$. Let $D$ be the largest amplitude of all the perturbation signals $d_i(t)$. Then it can be concluded that there exists a positive constant $\alpha_L$ such that:
\begin{align}
    \dot{L} \leq -\alpha_L L + \frac{\beta \|\bf{B}\| k_4 D^2}{2} .
\end{align}{}With $z = \left(\tilde{\bf{\eta}}, \tilde{\bf{\theta}}, \bf{x}, \tilde{\bf{u}}\right) \in \mathbb{R}^N \times \mathbb{R}^{m+N} \times \mathbb{R}^{n} \times \mathbb{R}^{m} $, let us define the set 
\begin{align} \Omega_{\gamma}=&\left\{z\ |\ L(z) \leq \gamma\right\}. \end{align}We choose $\gamma$ such that
\begin{align}
z \in \Omega_{\gamma} \Rightarrow \hat{\bf{\theta}} \in \Theta_1 \times \Theta_2 \times \dots \times \Theta_N.
\end{align}{}It then follows that the trajectories $\tilde{\bf{\eta}}, \tilde{\bf{\theta}}, \bf{x}, \tilde{\bf{u}} $ enter the set
\begin{align}
    \Omega_{\gamma_0}=\left\{z\ |\ L(z) \leq \frac{\beta \|\bf{B}\| k_4 D^2}{2\alpha_L} \right\}.
\end{align}
Therefore,  for $D$ chosen such that $\Omega_{\gamma_0} \subset \Omega_{\gamma}$, the set $\Omega_{\gamma_0}$, which is contained in a ball containing the point $(0, 0, \bf{x}^*, 0)$ with radius of order $\mathcal{O}\left(D^2\right)$, is exponentially stable for the closed-loop system . \hfill $\blacksquare$

\bibliographystyle{unsrt}
\bibliography{ifacconf} 

\begin{thebibliography}{24}
\providecommand{\natexlab}[1]{#1}
\providecommand{\url}[1]{\texttt{#1}}
\providecommand{\urlprefix}{URL }
\expandafter\ifx\csname urlstyle\endcsname\relax
  \providecommand{\doi}[1]{doi:\discretionary{}{}{}#1}\else
  \providecommand{\doi}{doi:\discretionary{}{}{}\begingroup
  \urlstyle{rm}\Url}\fi

\bibitem[{Adetola and Guay(2007)}]{adetola2007parameter}
Adetola, V. and Guay, M. (2007).
\newblock Parameter convergence in adaptive extremum-seeking control.
\newblock \emph{Automatica}, 43(1), 105--110.

\bibitem[{Adetola and Guay(2008)}]{adetola2008finite}
Adetola, V. and Guay, M. (2008).
\newblock Finite-time parameter estimation in adaptive control of nonlinear
  systems.
\newblock \emph{IEEE Transactions on Automatic Control}, 53(3), 807--811.

\bibitem[{Bauschke et~al.(2011)Bauschke, Combettes et~al.}]{bauschke2011convex}
Bauschke, H.H., Combettes, P.L., et~al. (2011).
\newblock \emph{Convex analysis and monotone operator theory in Hilbert
  spaces}, volume 408.
\newblock Springer, 2 edition.

\bibitem[{Bianchi and Grammatico(2019)}]{mattia2019dynamic}
Bianchi, M. and Grammatico, P. (2019).
\newblock A continuous-time distributed generalized nash equilibrium seeking
  algorithm over networks for double-integrator agents.
\newblock \emph{arXiv preprint arXiv:1910.11608}.

\bibitem[{De~Persis and Grammatico(2019)}]{de2019distributed}
De~Persis, C. and Grammatico, S. (2019).
\newblock Distributed averaging integral nash equilibrium seeking on networks.
\newblock \emph{Automatica}, 110, 108548.

\bibitem[{D{\"u}rr et~al.(2013)D{\"u}rr, Stankovi{\'c}, Ebenbauer, and
  Johansson}]{durr2013lie}
D{\"u}rr, H.B., Stankovi{\'c}, M.S., Ebenbauer, C., and Johansson, K.H. (2013).
\newblock Lie bracket approximation of extremum seeking systems.
\newblock \emph{Automatica}, 49(6), 1538--1552.

\bibitem[{Frihauf et~al.(2011)Frihauf, Krstić, and Basar}]{frihauf2011nash}
Frihauf, P., Krstić, M., and Basar, T. (2011).
\newblock Nash equilibrium seeking in noncooperative games.
\newblock \emph{IEEE Transactions on Automatic Control}, 57(5), 1192--1207.

\bibitem[{Gadjov and Pavel(2018)}]{gadjov2018passivity}
Gadjov, D. and Pavel, L. (2018).
\newblock A passivity-based approach to nash equilibrium seeking over networks.
\newblock \emph{IEEE Transactions on Automatic Control}, 64(3), 1077--1092.

\bibitem[{Ghaffari et~al.(2012)Ghaffari, Krsti{\'c}, and
  Ne{\v{s}}I{\'c}}]{ghaffari2012multivariable}
Ghaffari, A., Krsti{\'c}, M., and Ne{\v{s}}I{\'c}, D. (2012).
\newblock Multivariable newton-based extremum seeking.
\newblock \emph{Automatica}, 48(8), 1759--1767.

\bibitem[{Grammatico(2017)}]{grammatico2017dynamic}
Grammatico, S. (2017).
\newblock Dynamic control of agents playing aggregative games with coupling
  constraints.
\newblock \emph{IEEE Transactions on Automatic Control}, 62(9), 4537--4548.

\bibitem[{Guay and Dochain(2017{\natexlab{a}})}]{guay2017proportional}
Guay, M. and Dochain, D. (2017{\natexlab{a}}).
\newblock A proportional-integral extremum-seeking controller design technique.
\newblock \emph{Automatica}, 77, 61--67.

\bibitem[{Guay and Dochain(2017{\natexlab{b}})}]{guay}
Guay, M. and Dochain, D. (2017{\natexlab{b}}).
\newblock A proportional-integral extremum-seeking controller design technique.
\newblock \emph{Automatica}, 77, 61--67.

\bibitem[{Guay et~al.(2018)Guay, Vandermeulen, Dougherty, and
  McLellan}]{guay2018distributed}
Guay, M., Vandermeulen, I., Dougherty, S., and McLellan, P.J. (2018).
\newblock Distributed extremum-seeking control over networks of dynamically
  coupled unstable dynamic agents.
\newblock \emph{Automatica}, 93, 498--509.

\bibitem[{Krstić and Wang(2000)}]{krstic2000stability}
Krstić, M. and Wang, H.H. (2000).
\newblock Stability of extremum seeking feedback for general nonlinear dynamic
  systems.
\newblock \emph{Automatica}, 36(4), 595--601.

\bibitem[{Leblanc(1922)}]{leblanc}
Leblanc, M. (1922).
\newblock Sur l’electri"cation des chemins de fer au moyen de courants
  alternatifs de frequence elevee.
\newblock \emph{Revue Generale de l’Electricite}.

\bibitem[{Lin et~al.(2014)Lin, Qu, and Simaan}]{lin2014distributed}
Lin, W., Qu, Z., and Simaan, M.A. (2014).
\newblock Distributed game strategy design with application to multi-agent
  formation control.
\newblock In \emph{53rd IEEE Conference on Decision and Control}, 433--438.
  IEEE.

\bibitem[{Liu and Krsti{\'c}(2011)}]{liu2011stochastic}
Liu, S.J. and Krsti{\'c}, M. (2011).
\newblock Stochastic nash equilibrium seeking for games with general nonlinear
  payoffs.
\newblock \emph{SIAM Journal on Control and Optimization}, 49(4), 1659--1679.

\bibitem[{Nagurney and Zhang(2012)}]{nagurney2012projected}
Nagurney, A. and Zhang, D. (2012).
\newblock \emph{Projected dynamical systems and variational inequalities with
  applications}, volume~2.
\newblock Springer Science \& Business Media.

\bibitem[{Poveda and Teel(2017)}]{poveda2017framework}
Poveda, J.I. and Teel, A.R. (2017).
\newblock A framework for a class of hybrid extremum seeking controllers with
  dynamic inclusions.
\newblock \emph{Automatica}, 76, 113--126.

\bibitem[{Romano and Pavel(2019)}]{romano2019dynamic}
Romano, A.R. and Pavel, L. (2019).
\newblock Dynamic ne seeking for multi-integrator networked agents with
  disturbance rejection.
\newblock \emph{arXiv preprint arXiv:1903.02587}.

\bibitem[{Saad et~al.(2012)Saad, Han, Poor, and Basar}]{saad2012game}
Saad, W., Han, Z., Poor, H.V., and Basar, T. (2012).
\newblock Game-theoretic methods for the smart grid: Game-theoretic methods for
  the smart grid: An overview of microgrid systems, demand-side management, and
  smart grid communications.
\newblock \emph{IEEE Signal Processing Magazine}, 29, 86--105.

\bibitem[{Tan et~al.(2006)Tan, Ne{\v{s}}i{\'c}, and Mareels}]{tan2006non}
Tan, Y., Ne{\v{s}}i{\'c}, D., and Mareels, I. (2006).
\newblock On non-local stability properties of extremum seeking control.
\newblock \emph{Automatica}, 42(6), 889--903.

\bibitem[{Yi and Pavel(2019)}]{yi2019operator}
Yi, P. and Pavel, L. (2019).
\newblock An operator splitting approach for distributed generalized nash
  equilibria computation.
\newblock \emph{Automatica}, 102, 111--121.

\bibitem[{Yu et~al.(2017)Yu, Van Der~Schaar, and Sayed}]{yu2017distributed}
Yu, C.K., Van Der~Schaar, M., and Sayed, A.H. (2017).
\newblock Distributed learning for stochastic generalized nash equilibrium
  problems.
\newblock \emph{IEEE Transactions on Signal Processing}, 65(15), 3893--3908.

\end{thebibliography}

\end{document}